\documentclass[12pt]{amsart}
\usepackage{amsmath, amssymb, latexsym, amsthm}

\newcommand{\ed}{

\end{document}
}

\long\def\forget#1\forgotten{}
\newcommand{\eps}{\varepsilon}

\newenvironment{myrules}
{\begin{list}{}
{
 \setlength{\leftmargin}{0.8cm}
 \setlength{\labelwidth}{2cm}
 \setlength{\labelsep}{0.2cm}
 \setlength{\parsep}{0.5ex plus 0.2ex minus 0.1 ex}
 \setlength{\itemsep}{0.4ex plus 0.2 ex minus 0ex}
}}{\end{list}}
\newcommand{\comment}[1]{}

\renewcommand{\diamond}{\diamondsuit}
\newcommand{\next}[1]{{#1}^+}

\newcommand{\br}[1]{\bigl(#1\bigr)}

\newcommand{\Occ}{Oracle chain condition}
\newcommand{\occ}{oracle chain condition}
\newcommand{\Trap}{\op{Trap}}
\newcommand{\cc}{chain condition}
\newcommand{\fY}{\mathfrak{Y}}

\newcount\skewfactor
\def\mathunderaccent#1#2 {\let\theaccent#1\skewfactor#2
\mathpalette\putaccentunder}
\def\putaccentunder#1#2{\oalign{$#1#2$\crcr\hidewidth
\vbox to.2ex{\hbox{$#1\skew\skewfactor\theaccent{}$}\vss}\hidewidth}}
\def\name{\mathunderaccent\tilde-3 }

\renewcommand{\phi}{\varphi}
\newcommand{\fZ}{\mathfrak{Z}}
\newcommand{\fF}{\mathfrak{F}}
\newcommand{\fD}{\mathfrak{D}}

\newcommand{\Dfin}{\mathfrak{D}_\mathrm{fin}}

\renewcommand{\P}{\mathbb{P}}
\newcommand{\bbC}{\mathbb{C}}

\newcommand{\arx}[1]{\texttt{http://arxiv.org/abs/#1}}
\newcommand{\bq}{\begin{quote}}
\newcommand{\eq}{\end{quote}}

\newcommand{\inv}{^{-1}}
\newcommand{\Cantor}{{{}^\w 2}}

\newcommand{\op}{\operatorname}

\newcommand{\Null}{\mathcal{N}}

\newcommand{\M}{\mathcal{M}}

\newcommand{\ww}{{{}^{\w}\w}}
\newcommand{\roth}{[\w]^\w}

\newcommand{\Q}{\mathbb{Q}}
\newcommand{\R}{\mathbb{R}}

\newcommand{\cU}{\mathcal{U}}
\newcommand{\cG}{\mathcal{G}}
\newcommand{\Union}{\bigcup}

\renewcommand{\b}{\mathfrak{b}}
\renewcommand{\t}{\mathfrak{t}}

\renewcommand{\c}{\mathfrak{c}}
\renewcommand{\d}{\mathfrak{d}}

\newcommand{\g}{\mathfrak{g}}
\renewcommand{\i}{\item}

\renewcommand{\r}{\mathfrak{r}}
\renewcommand{\u}{\mathfrak{u}}
\newcommand{\h}{\mathfrak{h}}
\newcommand{\p}{\mathfrak{p}}

\newcommand{\s}{\mathfrak{s}}
\newcommand{\w}{\omega}

\newcommand{\nin}{\not\in}

\newcommand{\sbst}{\subseteq}

\newcommand{\sm}{\setminus}

\newcommand{\as}{\subseteq^*}
\renewcommand{\(}{\left(}
\renewcommand{\)}{\right)}
\newcommand{\<}{\langle}
\renewcommand{\>}{\rangle}

\newcommand{\dom}{\op{dom}}
\newcommand{\cov}{\mathsf{cov}}
\newcommand{\add}{\mathsf{add}}
\newcommand{\cof}{\mathsf{cof}}
\newcommand{\cf}{\mathsf{cf}}
\newcommand{\non}{\mathsf{non}}

\newtheorem{thm}{Theorem}[section]
\newtheorem{prop}[thm]{Proposition}

\newtheorem{lem}[thm]{Lemma}
\newtheorem{cor}[thm]{Corollary}

\theoremstyle{definition}
\newtheorem{choice}[thm]{Choice} 
\newtheorem{claim}[thm]{Claim}   
\newtheorem{defn}[thm]{Definition}
\newtheorem{notn}[thm]{Notation}
\theoremstyle{remark}

\newcommand{\be}{\begin{enumerate}}
\newcommand{\ee}{\end{enumerate}}
\newcommand{\bi}{\begin{itemize}}
\newcommand{\ei}{\end{itemize}}

\forget  
 \forgotten

\author{Heike Mildenberger}
\address{Universit\"at Wien,
Institut f\"ur Formale Logik, W\"ahringer Str.\ 25, 1090 Vienna,
Austria}
\email{heike@logic.univie.ac.at}

\author{Saharon Shelah}
\address{Einstein Institute of Mathematics, The Hebrew University of Jerusalem,
Givat Ram, 91904 Jerusalem, Israel, and
Mathematics Department, Rutgers University, 110 Frelinghuysen Road, NJ 08854-8019,
USA}
\email{shelah@math.huji.ac.il}

\author{Boaz Tsaban}
\address{Einstein Institute of Mathematics, The Hebrew University of Jerusalem,
Givat Ram, 91904  Jerusalem, Israel}
\email{tsaban@math.huji.ac.il}
\urladdr{http://www.cs.biu.ac.il/\~{}tsaban}


\thanks{The authors were partially supported by: The Austrian
``Fonds zur wissenschaftlichen F\"orderung'', grant no.\ 16334,
and the University of Helsinki (first author),
the Edmund Landau Center for Research in Mathematical Analysis and
Related Areas, sponsored by the Minerva Foundation, Germany (first and third author),
the United States-Israel Binational Science Foundation Grant no.\ 2002323
(second author), and the Golda Meir Fund (third author).
This is the second author's publication 847.}

\title[Covering and finite dominance]
{Covering the Baire space by families which are not finitely dominating}

\begin{document}
\begin{abstract}
It is consistent (relative to ZFC) that
each union of $\max\{\b,\g\}$ many families in the Baire space $\ww$
which are not finitely dominating is not dominating.
In particular, it is consistent that for each nonprincipal ultrafilter
$\cU$, the cofinality of the reduced ultrapower $\ww/\cU$ is greater than
$\max\{\b,\g\}$.
The model is constructed by \occ{} forcing,
to which we give a self-contained introduction.
\end{abstract}

\maketitle

\section{Introduction}
The undefined terminology used in this paper is as in \cite{Vaughan, BlassHBK}.
A family $Y\sbst\ww$ is \emph{finitely dominating} if
for each $g\in\ww$ there exist $k$ and $f_1,\dots,f_k\in Y$
such that $g(n)\le\max\{f_1(n),\dots,f_k(n)\}$ for all
but finitely many $n$.
The \emph{additivity number} for classes $\fY\sbst\fZ\sbst P(\ww)$ with $\Union\fY\nin\fZ$ is
$$\add(\fY,\fZ)=\min\{|\fF| : \fF\sbst\fY\mbox{ and }\Union\fF\nin\fZ\}.$$
Let $\fD$ (respectively, $\Dfin$) be the collection of all subsets of $\ww$ which are not
dominating (respectively, finitely dominating).
Define
$$\cov(\Dfin)=\min\{|\fF| : \fF\sbst\Dfin\mbox{ and }\Union\fF=\ww\}.$$
It is easy to see that $\add(\Dfin,\fD)=\cov(\Dfin)$, so we will
use this shorter notation.

In \cite{ShTb768} it is pointed out that
$$\max\{\b,\g\}\le\cov(\Dfin),$$
the inequality $\b\le\cov(\Dfin)$ being immediate from the
definitions, and the inequality $\g\le\cov(\Dfin)$ having been implicitly
proved in \cite[Theorem 2.2]{Mildenberger}.
(For the reader's convenience, we give a short proof for this
in Corollary \ref{gDfin}).
In \cite{ShTb768} it is shown that in all ``standard'' forcing extensions
(e.g., those appearing in \cite[\S 11]{BlassHBK}),
equality holds.
It is conjectured in \cite{ShTb768} that this equality
is not provable. We prove this conjecture.
In fact, we prove a stronger result:
Let $\M$ denote the ideal of meager sets
of real numbers.
\begin{thm}\label{main}
It is consistent (relative to ZFC) that $\aleph_1=\non(\M)=\g<\cov(\Dfin)=\cov(\M)=\c=\aleph_2$.
\end{thm}
The statement of Theorem \ref{main} determines the values of almost all
standard cardinal characteristics of the continuum in the model witnessing it:
If $\Null$ is the ideal of null sets of real numbers, then
by provable inequalities (see \cite{Vaughan, BlassHBK}), we have that
$\p,\t,\h,\add(\Null),\add(\M),\b,\s,\cov(\Null)$, and $\non(\M)$ are all equal to $\aleph_1$,
and
$\cov(\M),\non(\Null),\r,\d,\allowbreak\u,
\mathfrak{i},\cof(\M)$, and $\cof(\Null)$ are all equal to $\aleph_2$
in this model.

In \cite{ShTb768} it is shown that
for each nonprincipal ultrafilter $\cU$ on $\w$,
$\cov(\Dfin)\le\cof(\ww/\cU)$. 

\begin{cor}
It is consistent (relative to ZFC) that
for each nonprincipal ultrafilter $\cU$ on $\w$,
$\max\{\b,\g\}<\cof(\ww/\cU)$.
\end{cor}
This corollary partially extends the closely related
Theorems 3.1 and 3.2 of \cite{ShSt465},
which are proved using the same machinery:
\Occ{} forcing.

\section{Making $\cov(\Dfin)$ and $\cov(\M)$ large}

From now on, by \emph{ultrafilter} we always mean a nonprincipal
ultrafilter on $\w$. We will use the following convenient
characterization. For functions $f,g\in\ww$ and an ultrafilter
$\cU$ we write $f\le_\cU g$ for $\{n : f(n)\le g(n)\}\in\cU$.

\begin{lem}[\cite{ShTb768}]\label{filtercovDfin}
For each cardinal number $\kappa$, the following are equivalent:
\be
\i $\kappa<\cov(\Dfin)$;
\i For each $\kappa$-sequence $\<(\cU_\alpha,g_\alpha) : \alpha<\kappa\>$
with each $\cU_\alpha$ an
ultrafilter and each $g_\alpha\in\ww$
there exists $g\in\ww$ such that for each $\alpha<\kappa$,
$g_\alpha\le_{\cU_\alpha} g$.
\ee
\end{lem}

We first show how this characterization easily implies an assertion
made in the introduction.
\begin{defn}
For $A\in\roth$, define the function $\next{A}\in\ww$ by
$\next{A}(n)=\min\{k\in A: n<k\}$ for all $n$.
\end{defn}

\begin{cor}[\cite{Mildenberger}]\label{gDfin}
$\g\le\cov(\Dfin)$.
\end{cor}
\begin{proof}
We use Lemma \ref{filtercovDfin}.
Assume that $\kappa<\g$, and $(\cU_\alpha,g_\alpha)$, $\alpha<\kappa$, are given
with each $\cU_\alpha$ an ultrafilter and each $g_\alpha\in\ww$.
We must show that there exists $g\in\ww$ such that for each $\alpha<\kappa$,
$g_\alpha\le_{\cU_\alpha} g$.
We will use the following ``morphism''.

\begin{lem}\label{morph}
For each $f\in\ww$ and each ultrafilter $\cU$,
$$\cG_{\cU,f}=\{A\in\roth : f\le_\cU\next{A}\}$$
is groupwise dense.
\end{lem}
\begin{proof}
Clearly, $\cG_{\cU,f}$ is closed under taking almost subsets.
Assume that $\{[a_n,a_{n+1}): n\in\w\}$ is an interval partition of $\w$.
By merging consecutive intervals we may assume that for each $n$,
and each $k\in[a_n,a_{n+1})$, $f(k)\le a_{n+2}$.

Since $\cU$ is an ultrafilter, there exists $\ell\in\{0,1,2\}$ such that
$$A_\ell=\Union_n [a_{3n+\ell},a_{3n+\ell+1})\in\cU$$
Take $A=A_{\ell+2\bmod 3}$.
For each $k\in A_\ell$, let $n$ be such that $k\in[a_{3n+\ell},a_{3n+\ell+1})$.
Then $f(k)\le a_{3n+\ell+2}=\next{A}(k)$. Thus $A\in\cG_{\cU,f}$.
\end{proof}
Thus, we can take $A\in\bigcap_{\alpha<\kappa}\cG_{\cU_\alpha,g_\alpha}$
and $g=\next{A}$.
\end{proof}

How are we going force a large value for $\cov(\Dfin)$?
If $\cov(\Dfin)=\aleph_1$, then by Lemma \ref{filtercovDfin}
this is witnessed by a sequence
$\<(\cU_\alpha,g_\alpha) : \alpha<\aleph_1\>$.
To refute a single such witness, we will use the following forcing notion,
where $A_\alpha\in\cU_\alpha$ for each $\alpha<\aleph_1$.
\begin{defn}\label{Qdef}
Fix an ordinal $\gamma$.
Assume that $A_\alpha\in\roth$ and $g_\alpha\in\ww$ for $\alpha<\gamma$.
Define a forcing notion
$$\Q = \Q(A_\alpha,g_\alpha : \alpha<\gamma) = \{(n,h,F) : n\in\w, h\in{^n\w}, F\in[\gamma]^{<\aleph_0}\},$$
with $(n_1,h_1,F_1)\le (n_2,h_2,F_2)$ if
$n_1\le n_2$, $h_2\|n_1=h_1$, $F_1\sbst F_2$, and
$$\br{\forall \alpha\in F_1}\br{\forall n\in [n_1,n_2)\cap A_\alpha}\ g_\alpha(n)\le h_2(n).$$
\end{defn}
Observe that $\Q$ is $\sigma$-centered. 
$\Q$ is a restricted variant of the Hechler forcing.
Advanced readers are recommended to skip the proof of the
following lemma, which is the same as for the Hechler forcing.

\begin{lem}
Assume that $A_\alpha\in\roth\cap V$ and $g_\alpha\in\ww\cap V$ for each $\alpha<\gamma$.
Then for $\Q = \Q(A_\alpha,g_\alpha : \alpha<\gamma)$,
$V^\Q\models\br{\exists g\in\ww}\br{\forall\alpha<\gamma}\ A_\alpha\as \{n : g_\alpha(n)\le g(n)\}$.
\end{lem}
\begin{proof}
Assume that $G$ is a $\Q$-generic filter over $V$.
Let $g=\Union\pi_2[G]$, where $\pi_2$ denotes the projection on the second coordinate.
Clearly, $g$ is a partial function from $\w$ to $\w$.
By density arguments, we have that $g$ is as required.
To see this, consider first the sets
$$D_m = \{(n,h,F)\in\Q : m\le n\}$$
for $m\in\w$. Each $D_m$ is dense in
$\Q$: Assume that $(n,h,F)\in\Q$.
If $m\le n$ then $[n,m)=\emptyset$; therefore
$(n,h,F)\le(n,h,F\cup\{\alpha\})\in D_m$.
Otherwise, define $h':m\to\w$ by $h'(k)=h(k)$ for $k<n$, and
$h'(k)=\max\{f_\beta(k) : \beta\in F\}$ for $k\in [n,m)$.
Then $(m,h',F)$ is a member of $D_{m,\alpha}$ extending $(n,h,F)$.
The density of the sets $D_m$ implies that
$\dom(g)=\w$.
Moreover, for each $\alpha<\gamma$ the set
$$E_\alpha = \{(n,h,F)\in\Q : \alpha\in F\}$$
is dense in $\Q$ (for each condition $(n,h,F)$, $(n,h,F\cup\{\alpha\})$ is
a stronger condition which belongs to $E_\alpha$).
Now fix $\alpha<\gamma$ and choose an element $(n_0,h_0,F_0)\in G\cap E_\alpha$.
For each $n\in A_\alpha\sm n_0$ choose an element $(n_1,h_1,F_1)\in G\cap D_{n+1}$,
and a common extension $(n_2,h_2,F_2)$ of $(n_0,h_0,F_0)$ and $(n_1,h_1,F_1)$.
As $\alpha\in F_0$ and $n\in [n_0,n_2)\cap A_\alpha$, we have that
$g_\alpha(n)\le g(n)$. Since this holds for each $n\ge n_0$,
we have that $A_\alpha\as\{n : g_\alpha(n)\le g(n)\}$.
\end{proof}
Consequently, doing an iteration of forcing notions with the above forcing
used cofinally often, with $\gamma=\aleph_1$ and an appropriate book-keeping
will increase $\cov(\Dfin)$. We will be more precise in the proof of Theorem~\ref{2.9}.

Observe that the sets $A_\alpha$ played no special role and in fact we could take
$A_\alpha=\w$ for each $\alpha$ (in this case we obtain a dominating real).
However, this freedom to choose $A_\alpha$ will play a crucial role in the sequel,
where we would like to make sure that $\b$ (or $\non(\M)$) and $\g$ remain
small while we increase $\cov(\Dfin)$.

\medskip

We now make some easy observations concerning our planned forcing.
We will construct our model by a finite support iteration
$\<\P_\alpha,\Q_\alpha:\alpha<\aleph_2\>$ of c.c.c.\ forcing
notions $\Q_\alpha$ which add reals for cofinally many $\alpha<\aleph_2$.
Consequently, $V^\P$ satisfies $\c\ge\aleph_2$, where
$\P=\P_{\aleph_2}=\Union_{\alpha<\aleph_2}\P_\alpha$. The model
$V$ we begin with will satisfy $V=L$ (in fact,
$\diamond^*_{\aleph_1}$ and $\diamondsuit_{\aleph_2}(S^2_1)$, with
$S^2_1 = \{\alpha<\aleph_2 : \cf(\alpha) = \aleph_1\}$, are
enough). Consequently, $V$ satisfies $|\P|=\aleph_2=2^{\aleph_1}$.
Since $\P$ satisfies the c.c.c., (nice) $\P$-names for reals
are countable and therefore there are at most
$|\P|^{\aleph_0}=2^{\aleph_1}=\aleph_2$ names for reals in $\P$,
so $V^\P\models\c=\aleph_2$.

Since we are using a finite support iteration, Cohen reals are introduced cofinally often
along the iteration, and this is well known to imply $\cov(\M)\ge\aleph_2$ in the final model
(briefly: Each meager set in the final model is contained in an $F_\sigma$, thus Borel, meager set.
Each Borel set is coded by a real, and every real appears at a stage $\alpha<\aleph_2$, so Cohen reals
added later will not belong to the Borel meager set which is the interpretation of this code,
and since this property is absolute, they will not belong to the interpretation in the final model.
Since $\aleph_2$ is regular, the codes for $\aleph_1$ many Borel meager sets all appear at an intermediate
stage, so their union does not contain Cohen reals added later).
\begin{cor}
In the final model, $\cov(\M)=\c=\aleph_2$ holds.
\end{cor}

Now we show how to impose some more constraints on our iteration
$\< \P_\alpha, \Q_\alpha: \alpha<\aleph_2\>$
so that in $V^{\P_{\aleph_2}}$, $\cov(\Dfin) = \aleph_2$. Our
exposition follows closely the treatment of names given in
\cite{BurMil}.

\begin{choice}\label{2.8}
We fix a $\diamondsuit_{\aleph_2}(S_1^2)$-sequence $\< S_\delta : \delta \in S_1^2\>$
in the ground model.
The idea is that stationarily often $S_\delta$ will
guess a function
\begin{equation}\label{task}
f : (\aleph_1\times
\aleph_2)\cup\aleph_1 \to ([\aleph_2]^{\leq \aleph_0})^{\aleph_0}.
\end{equation}
(So for each $\delta<\aleph_2$ of cofinality $\aleph_1$,
$S_\delta :  ( \aleph_1\times \delta ) \cup \aleph_1 \to ([\delta]^{\leq\aleph_0})^{\aleph_0}$.)

We identify $\aleph_2$ with the partial order $\P_{\aleph_2}$ we are about to build.
Then $[\aleph_2]^{\leq\aleph_0}$ contains all of the maximal antichains. Thus
$([\aleph_2]^{\leq\aleph_0})^{\aleph_0}$ contains a name for each
subset of $\omega$ (which corresponds to an element of $\ww$).
Now any sequence
$$\< (\cU_\alpha, g_\alpha) : \alpha<\aleph_1\>$$
in the extension has a ground model function $f: ( \aleph_1\times
\aleph_2)\cup \aleph_1 \to ([\aleph_2]^{\leq\aleph_0})^{\aleph_0}$,
such that $f(\alpha)$ is a name for $g_\alpha$ and $f(\alpha,\cdot)$ is a name for an enumeration
of the elements of
$\cU_\alpha$.

For each $f$ as in Equation~\eqref{task},
$$\{\delta \in S_1^2 : S_\delta= f\restriction \delta\}$$
is stationary in $\aleph_2$.
We will inductively define an $\aleph_2$-stage finite support iteration
and an injection function $F_\delta : \P_\delta \to \aleph_2$ for $\delta < \aleph_2$
such that the range of each $F_\delta$ is an initial segment of $\aleph_2$ which includes
$\delta$, and for $\eps<\delta <\aleph_2$, $F_\eps \subseteq F_\delta$.

For $\delta <\aleph_2$ we will denote by  ${\rm name}(S_\delta)$ the sequence of
$\aleph_1$ sets of reals $\cU_\alpha$ and of $\aleph_1$ reals $g_\alpha$
of the form
\begin{multline*}
\< (\{\bigcup_{n\in\omega} \{n\}
\times F_\delta^{-1}(S_\delta(\alpha,\xi)(n)) : \xi <\delta \},
\bigcup_{n\in\omega}\{n\} \times F_\delta^{-1}(S_\delta(\alpha)(n)))
:\\ \alpha< \aleph_1\>.
\end{multline*}

At stage $\delta \in S_1^2$ in the construction, if
$\Vdash_{\P_{\delta}}\mbox{``} {\rm name}(S_\delta)$ is a sequence of
$\aleph_1$ ultrafilters and $\aleph_1$ functions'', then we can take $\P_\delta$-names
$A_\alpha$,
$\alpha<\aleph_1$, such that
$\Vdash_{\P_\delta} A_\alpha \in (\cU_\alpha)\restriction \delta$,
which means
$\Vdash_{\P_\delta} \mbox{``}A_\alpha$ is in the first component of ${\rm name}(S_\delta)$''.
\end{choice}

\begin{thm}\label{2.9}
Let $V\models \diamond_{\aleph_2}(S_1^2)$ and
let $\P_{\aleph_2}$ be any forcing as in Choice~\ref{2.8}.
 Then $V^{\P_{\aleph_2}} \models \cov(\Dfin) = \aleph_2$.
\end{thm}

\begin{proof} If $\Vdash_{\P_{\aleph_2}} \mbox{``} \< (\cU_\alpha,g_\alpha)
: \alpha<\aleph_1\>$ is a sequence of functions and ultrafilters'',
then at club many stages $\delta$ the restriction of the names to
$\delta$ is also forced
to be a sequence of ultrafilters in $V^{\P_{\delta}}$.
For a proof of this (even in the countable support proper scenario) see
\cite{BlSh}.
But the restriction of the name to $\delta$ is guessed by ${\rm name}(S_\delta)$
for stationarily many $\delta$'s in this club.
So at such a stage $\delta$ the forcing $\Q_\delta$ adds a function $h$ such that
$g_\alpha \leq _{\cU_\alpha} h$ for all $\alpha<\aleph_1$ and this shows that
the sequence was not a witness for $\cov(\Dfin)= \aleph_1$.
\end{proof}

\section{Interlude: \Occ{} forcing}

Usually, the major difficulty in forcing inequalities between
combinatorial cardinal characteristics of the continuum is to make
sure that those which are required to be smaller ($\non(\M)$ and $\g$ in our case)
indeed remain small in the generic extension.
In this section we describe one such method, which is suitable for our purposes:
\Occ{} forcing \cite[Chapter IV]{Sh:f} (see also \cite{Burke, BurMil}).

\Occ{} forcing is a method for forcing with $\aleph_2$-stage finite
support iteration,
in such a way that some prescribed intersections of $\aleph_1$ many
(descriptively nice) sets which are empty in an intermediate model
remain empty in the final model.

\begin{defn}
An \emph{oracle} (or \emph{$\aleph_1$-oracle})
is a sequence $\bar{M} = \<M_\delta:\delta\mbox{ limit}<\aleph_1\>$
of countable transitive models of a sufficiently large finite portion of ZFC (henceforth denoted ZFC$^*$),
such that for each $\delta$, $\delta\in M_\delta$ is countable in $M_\delta$,
and for each $A\sbst\aleph_1$, the set
$$\Trap_{\bar M}(A)=\{\delta<\aleph_1 : \delta\mbox{ is a limit ordinal, and }A \cap \delta \in M_\delta \}$$
is a stationary subset of $\aleph_1$.
\end{defn}
Clearly, $\diamond$ implies the existence of an oracle.
The sets $\Trap_{\bar M}(A)$ generate a filter $\Trap_{\bar M}$,
which is normal and proper.
Moreover, for each $A,B\sbst\aleph_1$,
there exists $C\sbst\aleph_1$ such that $\Trap_{\bar M}(C)=\Trap_{\bar M}(A)\cap\Trap_{\bar M}(B)$.

\begin{notn}
Assume that $\P\sbst\Q$ are forcing notions, and $N$ is a set.
Then $\P<_N\Q$ means: Every predense subset of $\P$ which belongs to
$N$ is predense in $\Q$.
\end{notn}

\begin{lem}\label{basicpredense}
~\be
\i $<_N$ is transitive,
\i If $N\sbst N'$, then $\P<_{N'}\Q$ implies $\P<_N\Q$;
\i If $\Q=\Union_{\alpha<\beta}\Q^\alpha$ and
$\P<_N\Q^\alpha$ for each $\alpha$, then
$\P<_N\Q$.\hfill\qed
\ee
\end{lem}

\begin{defn}\label{occdef}
Assume that $\bar M$ is an oracle.
A forcing notion $\P$ satisfies the \emph{$\bar M$-\cc{}}
if there exists an injection $\iota:\P\to\aleph_1$, such that
$$\{\delta<\aleph_1 : \delta\mbox{ is a limit ordinal, and }\iota^{-1}[\delta]<_{M_{\delta,\iota}}\P\}
\in\Trap_{\bar M},$$
where $M_{\delta,\iota}=\{\iota^{-1}[A] : A\sbst\delta\mbox{ and }A\in M_\delta\}$.
\end{defn}

Thus each countable forcing notion satisfies the $\bar M$-\cc{}, and if $\P$
satisfies the $\bar M$-\cc{}, then $\P$ has the c.c.c., and $|\P|\le\aleph_1$.
The definition of the $\bar M$-\cc{} can be extended to forcing notions
of cardinality $\aleph_2$ \cite[IV.1.5]{Sh:f}; however this is not needed here.

Proving the $\bar M$-\cc{} according to Definition \ref{occdef} is rather inconvenient.
We give a useful method to verify the $\bar M$-\cc{}.

\begin{prop}\label{occrecipe}
Assume that $\bar M$ is an oracle,
$\P=\Union_{\alpha<\aleph_1}\P^\alpha$,
for each $\alpha<\aleph_1$, $\iota_\alpha$ is a bijection from $\P^\alpha$ onto a countable ordinal,
and $\<N_\alpha : \alpha<\aleph_1\>$ is a sequence of countable transitive models of ZFC$^*$,
such that the following conditions hold:
\be
\i For each $\alpha<\beta<\aleph_1$,
\be
\i $\P^\alpha\sbst\P^\beta$ with $\P^\beta\sm\P^\alpha$ countably infinite,
\i $\iota_\alpha\sbst \iota_\beta$; and
\i $N_\alpha\sbst N_\beta$. 
\ee
\i For each (large enough) $\alpha<\aleph_1$,
\be
\i $\iota_\alpha:\P^\alpha\to\w\alpha$ is bijective,
\i $M_{\omega\alpha}, \<\P^\alpha,\le_{\P^\alpha}\>, \iota_\alpha\in N_{\alpha}$; and
\i $\P^\alpha <_{N_\alpha}\P^{\alpha+1}$.
\ee
\ee
Then $\P$ satisfies the $\bar M$-\cc{}.
\end{prop}
\begin{proof}
Using Lemma \ref{basicpredense}, we get
by induction on $\beta$ that
for each $\alpha\le\beta\le\aleph_1$, $\P^\alpha <_{N_\alpha}\P^{\beta}$.
In particular, $\P^\alpha <_{N_\alpha}\P$ for each $\alpha$.
Define $\iota=\Union_{\alpha<\aleph_1}\iota_\alpha$.
Then $\iota:\P\to\aleph_1$ is an injection.

Assume that $\delta<\aleph_1$ is a (large enough)
limit ordinal, and let $\alpha$ be such that $\delta=\w\alpha$.
Then
$$\iota\inv[\delta] = \iota\inv[\w\alpha] = \iota_\alpha\inv[\w\alpha]=\P^\alpha.$$
Assume that $A\sbst\delta$, $A\in M_\delta$, and
$\iota\inv[A]=\iota_\alpha\inv[A]$ is predense in $\P^\alpha$.
As $\iota_\alpha\in N_\alpha$, $\iota_\alpha\inv[A]\in N_\alpha$.
As $\P^\alpha <_{N_\alpha}\P$, $\iota_\alpha\inv[A]$ is predense in $\P$.

This shows that for \emph{all} (large enough) limit ordinals $\delta<\aleph_1$,
$\iota\inv[\delta]<_{M_{\delta,\iota}}\P$. Obviously, this implies the
requirement in Definition \ref{occdef}.
\end{proof}

Proposition \ref{occrecipe} gives us a recipe for verifying the $\bar M$-\cc{}:
Construct $\P$ by inductively constructing $\P^\beta$,
such that (1)(a) holds. If $\beta$ is a limit, take $\P^\beta=\Union_{\alpha<\beta}\P^\alpha$.
Otherwise $\beta=\alpha+1$ and $\P^\alpha$ is defined.
Then there exists $\iota_\beta$ such that (1)(b) and (2)(a)
hold. Choose $N_\alpha$ as in (1)(c) and (2)(b) (and containing some other elements if needed),
and use $N_\alpha$ to define $\P^{\alpha+1}$ such that (2)(c) holds (this is the only tricky
part in the construction). We can simplify the last step in this recipe a bit further.

\begin{lem}
Assume that $N$ is a transitive model of ZFC$^*$,
such that $\<\P,\le_\P\>\in N$. Then:
$\P<_N\Q$ if, and only if, each open dense subset of $\P$
which belongs to $N$ is predense in $\Q$.
\end{lem}
\begin{proof}
We need to prove $(\Leftarrow$).
Assume that $I\in N$ is predense in $\P$.
Then $I^* = \{p\in\P : \(\exists q\in I\) p\ge q\}\in N$,
and is open and dense in $\P$.
Thus, $I^*$ is predense in $\Q$, and therefore
$I$ is predense in $\Q$ as well.
\end{proof}

\begin{cor}\label{easier}
(2)(c) in Proposition \ref{occrecipe} can be replaced by:
\bi
\i[(2)(c$'$)] Each open dense subset of $\P^\alpha$
which belongs to $N_\alpha$ is predense in $\P^{\alpha+1}$.
\ei
\end{cor}

The following theorem exhibits the importance of the \occ{} for a single step forcing.

\begin{thm}[{\cite[IV.2.1]{Sh:f}}]\label{omit}
Assume that $V\models\diamond$, and
$\phi_\alpha(x)$, $\alpha<\aleph_1$, are $\Pi^1_2$ formulas\footnote{That is,
formulas of the form $\(\forall a\in\R\)\(\exists b\in\R\)\ \psi$,
where $\psi\in L_{\aleph_1,\aleph_0}$
($L_{\aleph_1,\aleph_0}$ is the extension of the first order
language by allowing countable conjunctions).}
(possibly with real parameters),
and $$V\models\lnot\(\exists x\)\(\forall\alpha<\aleph_1\)\ \phi_\alpha(x).$$
If this continues to hold when we add a Cohen real to $V$,
then
 there exists an oracle $\bar M$ such that for each forcing notion $\P$
satisfying the $\bar M$-\cc{},
$V^\P\models\lnot\(\exists x\)\(\forall\alpha<\aleph_1\)\ \phi_\alpha(x)$.
\end{thm}

The following consequence can be derived from Theorem \ref{omit}.

\begin{lem}[{\cite[IV.2.2]{Sh:f}}]\label{nonmgroracle}
Assume that $\diamond$ holds in $V$. There is an oracle $\bar M$ in $V$
such that for each $\P$ satisfying the $\bar M$-c.c.,
if, in $V$, $A$ is a nonmeager set of reals, then
$A$ is nonmeager in $V^\P$. Consequently, $V^\P\models\non(\M)=\aleph_1$.
\end{lem}

\Occ{} can (and is intended to) be used with finite support iterations.

\begin{lem}[{\cite[IV:3.2--3.3]{Sh:f}}]\label{iterate}
Assume that $\bar M$ is an oracle.
\be
\i\label{eachP}
For a finite support iteration $\<\P_\alpha,\name{\Q}_\alpha:\alpha<\gamma\>$,
if each $\P_\alpha$ satisfies the $\bar M$-\cc{}, then so does $\P_\gamma=\Union_{\alpha<\gamma}\P_\alpha$.
\i\label{eachQ}
If $|\P|=\aleph_1$, and $\P$ satisfies the $\bar M$-\cc{} (in $V$),
then in $V^\P$ there is an oracle $\bar M^*$ such that for each $\Q\in V^\P$ satisfying the
$\bar M^*$-\cc{}, $\P\star\name{\Q}$ satisfies the $\bar M$-\cc{} (in $V$).
\ee
\end{lem}

Consider a finite support iteration $\<\P_\alpha,\name{\Q}_\alpha : \alpha<\aleph_2\>$ of forcing
notions, and let $\P=\Union_{\alpha<\aleph_2}\P_\alpha$.
Assume that we wish to use Theorem \ref{omit} for $\P$.
Then by Lemma \ref{iterate}(\ref{eachP}), it suffices to make sure that
each $\P_\alpha$ satisfies the $\bar M$-\cc{}.
By Lemma \ref{iterate}(\ref{eachQ}), this amounts to
choosing each $\Q_\alpha$ in a way that it satisfies the \occ{} for the oracle
$\bar M^*$ corresponding to the oracle $\bar M$ given in Theorem \ref{omit} for $\P_\alpha$.

The nice thing is that we need not worry what exactly are these oracles,
as long as we can make sure that for any prescribed oracle $\bar M$,
the forcing notion $\Q_\alpha$ used in the iteration
can be chosen so that it satisfies the $\bar M$-\cc{}.

We sometimes have to make more than one oracle commitment.
In fact, we may wish to add new commitments cofinally often along the iteration
(indeed, we do that in the proof of Theorem~\ref{stronger}).
This can be achieved by coding all of the oracles of interest (those introduced
in earlier stages of the iteration as well as the new ones required in the
current iteration) in a single oracle.
Since the length of the iteration is $\aleph_2$, the following lemma tells that
this is possible.
\begin{lem}[{\cite[IV.3.1]{Sh:f}}]
If $\bar M_\alpha$, $\alpha<\aleph_1$, are oracles in $V$, then
there exists a single oracle $\bar M$ such that for each $\P$ satisfying
the $\bar M$-\cc{}, $\P$ satisfies the $\bar M_\alpha$-\cc{} for each $\alpha$.
\end{lem}

\section{Keeping $\non(\M)$ small}    

The main lemma needed to carry out our constructions is the following.
\begin{lem}\label{mainlemma}
Assume that $\bar M$ is an oracle, and
for each $\alpha<\aleph_1$, $\cU_\alpha$ is an ultrafilter
and $g_\alpha\in\ww$. Then there exist sets $A_\alpha\in\cU_\alpha$, $\alpha<\aleph_1$,
such that $\Q=\Q(A_\alpha,g_\alpha : \alpha<\aleph_1)$ (Definition \ref{Qdef})
satisfies the $\bar M$-\cc{}.
\end{lem}
\begin{proof}
We use Proposition \ref{occrecipe} and the remarks following it
(with $\P$ replaced by $\Q$ everywhere).
We choose $A_\alpha$ by induction on $\alpha$.
At stage $\alpha$ we define
$$\Q^\alpha = \Q(A_\beta,g_\beta : \beta<\alpha)$$
(so at the end, $\Q=\Union_{\alpha<\aleph_1}\Q^\alpha$ and (1)(a) is guaranteed)
and $\iota_\alpha$ as in (1)(b) and (2)(a),
then we choose $N_\alpha$ such that $N_\beta\sbst N_\alpha$ for each $\beta<\alpha$,
and $g_\alpha\in N_\alpha$ and (2)(b) holds.

Recall that $N_\alpha$ is countable, so we can choose an increasing sequence $\<a_k:k\in\w\>$
of natural numbers such that for each $g\in N_\alpha$, $g(a_k)<a_{k+1}$ for all but
finitely many $k$ (to obtain such a sequence, take an increasing function $f\in\ww$ which
dominates all members of $\ww\cap N_\alpha$, and define $a_k=f^k(0)$).
Since $\cU_\alpha$ is an ultrafilter, there exists $\ell\in\{0,1\}$ such that
$$A_\alpha := \Union_{k\in\w}[a_{2k+\ell},a_{2k+1+\ell})\in\cU_\alpha.$$
It remains to show that this definition guarantees (2)(c), that is,
$\Q^\alpha<_{N_\alpha}\Q^{\alpha+1}$.
We will use Corollary \ref{easier} for that.
Assume that $D\in N_\alpha$ is an open dense subset of $\Q^\alpha$,
and $p=(n,h,F)\in\Q^{\alpha+1}\sm\Q^\alpha$ (so $\alpha\in F$).
Define, for each $m>n$, $h_m:m\to\w$ by
$$h_m(k) = \begin{cases}
h(k) & k<n\\
\max\{g_\beta(k) : \beta\in F\} & n\le k
\end{cases}$$
Then $(n,h,F)\le (m,h_m,F)$,
and in particular $(n,h,F\sm\{\alpha\})\le (m,h_m,\allowbreak F\sm\{\alpha\})$.
Note that the mapping $m\mapsto h_m$ belongs to $N_\alpha$.

Define $f:\w\to\w$ by letting $f(k)$ be the minimal
$m$ such that there exists an element $(m,\tilde h,\tilde F)\in D$
which extends $(k,h_k,F\sm\{\alpha\})$.
Then $f\in N_\alpha$, so there exists
$k$ such that $m:=f(a_{2k+\ell-1})<a_{2k+\ell}$.
Let $q_0=(a_{2k+\ell-1},h_{a_{2k+\ell-1}},F\sm\{\alpha\})$.
By the definition of $f$, there exists
$q_1:=(m,\tilde h,\tilde F)\in D$ which extends
$q_0$.
Let $q_2=(m,\tilde h,\tilde F\cup\{\alpha\})\in\Q^{\alpha+1}$.

Then $q_1\le q_2$ since they share the same domain.
Since $q_1\in D$, it remains to show that $(n,h,F)\le q_2$.
$(n,h,F\sm\{\alpha\})\le q_0\le q_1$; thus
$(n,h,F\sm\{\alpha\})\le q_2$, and hence it suffices to show that
for each $i\in[n,m)\cap A_\alpha$, $g_\alpha(i)\le\tilde h(i)$.
But since $A_\alpha\cap[a_{2k+\ell-1},a_{2k+\ell})=\emptyset$,
$[n,m)\cap A_\alpha\sbst [n,a_{2k+\ell-1})$, and
if $i\in[n,a_{2k+\ell-1})$, then
$\tilde h(i)=h_{a_{2k+\ell-1}}(i)=\max\{g_\beta(i) : \beta\in F\}\ge g_\alpha(i)$,
since $\alpha\in F$, and we are done.
\end{proof}

By Lemma \ref{iterate}, Lemma \ref{mainlemma} will enable us to keep $\non(\M)$ small.
We now turn to the problem of keeping $\g$ small.

\section{Keeping $\g$ small}

First we state a sufficient condition for $\mathfrak g$ being small.

\begin{lem}\label{gcond} 
Assume that $\{Y_\zeta : \zeta < \c\}\sbst\roth$, and $\kappa$ is a cardinal such that:
\be
\i For each meager set ${\bf B}\sbst\roth$, $|\{ \zeta : Y_\zeta \not\in {\bf B}\}| = \c$.
\i For each $B\in\roth$, $|\{\zeta < \c : B \as Y_\zeta\}| <\kappa$.
\ee
Then $\g\leq \kappa$.
\end{lem}
\begin{proof}
By a result of Blass \cite{BlassHBK}, $\g\leq \cf(\c)$, so
we can assume that $\kappa \leq \cf(\c)$.
We now define $\kappa$ sets and then show that they are groupwise dense and that their
intersection is empty.

Let $\< \bar{n}^\zeta : \zeta < \c\>$
list all strictly increasing sequences of natural numbers, each sequence appearing
cofinally often.
By induction on $\zeta < \c$
we choose $\eps_\zeta \leq \kappa$, $\gamma_\zeta < \c$ and  $C_\zeta \in\roth$
as follows.

If there is some $\eps < \kappa$ such that for each $\xi < \zeta$
with $\eps_\xi = \eps$ we have $[n^\zeta_i,n^\zeta_{i+1}) \not\subseteq C_\xi$
for all but finitely many $i$, then we take as $\eps_\zeta$ the minimal such $\eps$.
By the assumption (1), we can choose $\gamma_\zeta$ to be the minimal $\gamma <\c$
such that $\gamma \neq \gamma_\xi$ for all $\xi <\zeta$ and there are infinitely many $i$
such that $[n^\zeta_i,n^\zeta_{i+1})\subseteq Y_\gamma$.
In this case we set $C_\zeta = \bigcup\{ [n^\zeta_i,n^\zeta_{i+1}) : i \in \omega, [n^\zeta_i,n^\zeta_{i+1})
\subseteq Y_{\gamma_\zeta}\}$.
Otherwise we set $\eps_\zeta = \kappa$ and $C_\zeta = \omega$.

For each $\xi <\kappa$, define
$$\cG_\xi = \{ B \in\roth : (\exists \zeta < \c)\ \eps_\zeta \geq \xi\mbox{ and } B \as C_\zeta\}.$$
We show that each $\cG_\xi$ is groupwise dense.
Clearly, it is closed under almost subsets.
Let an increasing sequence $\bar{n}$ be given.
Then for each $\nu < \xi$, there is by our construction some
$\zeta(\nu) < \c$ such that $\eps_{\zeta(\nu)} = \nu$ and $[n_i,n_{i+1})\subseteq C_{\zeta(\nu)}$
for infinitely many $i$.
As $\kappa \leq\cf(\c)$, $\zeta(*)=\sup\{ \zeta(\nu) : \nu < \xi \} <\c$.
By the choice of $\< \bar{n}^\zeta : \zeta < \c \>$ there is some
$\beta \in (\zeta(*),\c)$ such that $\bar{n}^\beta = \bar{n}$. So $\eps_\beta \geq \xi$,
and
$\bigcup \{ [n^\beta_i,n^\beta_{i+1}) :
[n^\beta_i,n^\beta_{i+1}) \subseteq Y_{\gamma_\beta} \} = C_\beta \in\cG_\xi$.

To see that $\bigcap\{ \cG_\xi : \xi<\kappa \}= \emptyset$, assume
that $B$ is infinite and for each $\xi$, $B \in \cG_\xi$.
Then for each $\xi <\kappa$, there is $\beta_\xi < \c$ such that
$\eps_{\beta_\xi}= \xi$ and $B \as C_{\beta_\xi} \subseteq Y_{\gamma_{\beta_\xi}}$.
Since $\kappa $ is regular, we can thin out and assume that
if  $\xi_1 <\xi_2$, then $\eps_{\beta_{\xi_1}} \neq \eps_{\beta_{\xi_2}}$.
Thus we have that for $\xi_1 < \xi_2$, $\beta_{\xi_1} \neq
\beta_{\xi_2}$, and hence $\gamma_{\beta_{\xi_1}} \neq \gamma_{\beta_{\xi_2}}$.
Consequently, $|\{ \gamma_{\beta_\xi} : \xi < \kappa \}| =\kappa$.
But $\{ \gamma_{\beta_\xi} : \xi < \kappa \} \subseteq \{\zeta < \c : B \as Y_\zeta\}$,
contradicting the assumption (2).
\end{proof}

As we already stated in the previous sections, we shall use a finite support iteration
$\<\P_\delta,\Q_\delta, : \delta<\aleph_2\>$ of c.c.c.\ forcing notions,
and choose constant or increasing oracles $\bar{M}^\delta$, such that
$\P_\delta$ has the $\bar{M}^\delta$-\cc{} for each $\delta$.
We start with a ground model satisfying
$\diamondsuit^*_{\aleph_1}$ and $\diamondsuit_{\aleph_2}(S^2_1)$.
Let $\< S_\delta : \delta \in S_1^2\>$ be a $\diamondsuit_{\aleph_2}(S_1^2)$-sequence.

There are three possibilities for $\Q_\delta$.
If $\cf(\delta) = \aleph_0$ or if $\delta$ is a successor,
then $\Q_\delta$ is the Cohen forcing.

If $\cf(\delta)= \aleph_1$ and $\Vdash_{{\P_\delta}}\mbox{``}{\rm name}(S_\delta)$
is a sequence of
ultrafilters $\cU_\alpha$ and of  functions $g_\alpha$,
$\alpha < \aleph_1$'', then
we choose $A_\alpha$, $\alpha < \aleph_1$ as in Lemma \ref{mainlemma}
but with additional provisos
and force with $\Q_\delta = \Q(\< A_\alpha ,g_\alpha : \alpha < \aleph_1 \>)$.
For the premise of this sentence we shortly say:
$S_\delta$ guesses $\< (\cU_\alpha, g_\alpha) : \alpha < \aleph_1\>$.
Otherwise, we set $\Q_\delta = \{0\}$.

\begin{defn}\label{5.2}
For $\gamma \leq \aleph_2$ we consider the class
${\mathcal K}_\gamma$  of \emph{$\gamma$-approximations}
$$\< (\P_\delta, \name{\Q_\delta}, \bar{M}^\delta,W_1, W_2) : \delta < \gamma\>$$
with the following properties:
\begin{myrules}
\item[(a)] $\< \P_\delta, \name{\Q_\delta} : \delta < \gamma\>$ is
a finite support iteration of partial orders such that for each $\delta < \gamma$,
$|\P_\delta| \leq \aleph_1$.
\item[(b)] $\< \bar{M}^\delta : \delta < \gamma \>$
is a constant sequence of oracles such that for all
$\delta$, $\P_\delta$ satisfies the $\bar{M}^\delta$-\cc{} and for $\delta+1 <\gamma$,
$\Vdash_{\P_\delta} \mbox{``}\name{\Q_\delta}$
satisfies the $(\bar{M}^{\delta+1})^*$-c.c.''
(as in Lemma \ref{iterate}(2)).
The constant value of the oracle sequence is
some oracle $\bar{M}$ as in Lemma~3.9, keeping $\cov(\M)= \aleph_1$.
\item[(c)] $W_1, W_2 \subseteq \aleph_2 \setminus S_1^2$, $W_1$ and  $W_2$ are
disjoint and if $\gamma$ is a limit of cofinality $\aleph_1$, then
$W_1 \cap \gamma$, $W_2\cap \gamma$ are both  cofinal in $\gamma$.
\item[(d)] If $\beta \in (W_1 \cup W_2)\cap \gamma$ then $\name{\Q_\beta}$ is the Cohen forcing
adding the real $\name{r_\beta} \in \Cantor$.
\item[(e)] If $\delta \in S_1^2\cap \gamma$ and $S_\delta$ guesses $\< (\cU_\alpha(\delta),
g_\alpha (\delta))
: \alpha <\aleph_1 \>$, then there is some strictly increasing  enumeration
 $\< \zeta_\alpha(\delta) :  \alpha < \aleph_1\>$ of a cofinal part of
$W_2 \cap \delta$, and  for every  $\alpha < \aleph_1$ there is $\ell_{\zeta_\alpha(\delta)}
 \in \{0,1\}$  such that
 $Y_{\zeta_\alpha(\delta)}^{\ell_{\zeta_\alpha(\delta)}}:=
r_{\zeta_\alpha(\delta)}^{-1}(\{\ell_{\zeta_\alpha(\delta)}\})
 \in \mathcal U_\alpha$, and
$\Q_\delta = \Q(Y_{\zeta_\alpha(\delta)}^{\ell_{\zeta_\alpha(\delta)}},
 g_\alpha(\delta) : \alpha < \aleph_1)$.\footnote{The $\zeta_\alpha(\delta)$,
$\alpha <\aleph_1$,
chosen here do not have to be coherent when regarding different $\delta$'s
and  we index them with $\delta$ because we need it.
Strictly speaking the $\ell_{\zeta_\alpha(\delta)}$ is a function
$\ell_{\zeta_\alpha(\delta)}(\delta)$. And also strictly speaking we
should index by $\gamma$ as well, but we are suppressing this because we are anyway
only working with end extensions when increasing $\gamma$.}
 \item[(f)] For all $\delta\leq\gamma$,
 $\Vdash_{\P_\delta} \mbox{``}(\forall A \in \roth)\ \{ \beta \in W_1 \cap \delta
: A \as \name{Y_\beta^1} \}$ is at most countable.''\footnote{Here it is $W_1$.
We use the Cohens in $W_2$ to build the forcings of type
$\Q_\delta =\Q(Y_{\zeta_\alpha(\delta)}^{\ell_{\zeta_\alpha}(\delta)},
g_\alpha(\delta) : \alpha <\aleph_1)$ and the Cohens $Y_\zeta^1$,
 $\zeta \in W_1$, to build the  $Y_\zeta$'s as in Lemma~\ref{gcond}.}
Here, for $\delta = \gamma$ limit, $\P_\gamma$ is the direct limit of $\langle
\P_\beta: \beta<\gamma\rangle$, and for $\delta = \gamma= \beta+1$,
$\P_\gamma = P_\beta \star \name{{\mathbb Q}}_\beta$.
\end{myrules}
\end{defn}

With the help of several lemmas we will prove the following.

\begin{thm}\label{goodclass}
If $V \models \diamondsuit_{\aleph_1}^*\mbox{ and }\diamondsuit_{\aleph_2}(S_1^2)$,
then for each $\gamma \leq \aleph_2$,  ${\mathcal K}_\gamma$ is not empty.
\end{thm}

Let $V$ fulfill the premises and let
$\P_{\aleph_2}$ be the direct limit of the first components of an $\aleph_2$-approximation.
If $G$ is a $\P_{\aleph_2}$-generic filter and
$\name{Y_\zeta^1}[G_{\aleph_2}] = Y_\zeta$ for $\zeta \in W_1$, then
we have in the final model a sequence $\< Y_\zeta : \zeta < \c\>$
as in Lemma~\ref{gcond} with $\kappa=\aleph_1$.

\begin{cor}
$V^{\P_{\aleph_2}} \models \cov(\M)= \g= \aleph_1<\cov(\Dfin)=\aleph_2$.
\end{cor}

We prove Theorem~\ref{goodclass} by induction on $\gamma$
and we shall work with end extensions.
For some $\gamma$'s, one has to work to show item (e). We will do this in our
first lemma.
For all $\gamma$'s but maybe the successor steps of points not in $S_1^2$, one has to work
to show that item (f) can be preserved in the induction. This will be
done in the last three lemmas.

\begin{lem}\label{5.5}
Consider a successor $\gamma=\delta+1$, $\delta \in S_1^2$. Given any $\aleph_1$-oracle
$(\bar{M}^{\delta+1})^*$,  the  sequence $\< \zeta_\alpha(\delta) : \alpha <\aleph_1\>$
can be chosen as in (e) so that the forcings given in item (e) have
the $(\bar{M}^{\delta+1})^*$-c.c.
\end{lem}

\begin{proof}
This is a variation of Lemma~4.1. We suppress some of the $\delta$'s.
We choose $\< \zeta_\alpha :  \alpha < \aleph_1\>$ enumerating
$W_2 \cap \delta$  so that, given the oracle
$(\bar{M}^{\delta+1})^*= \< N_\alpha : \alpha < \aleph_1 \>$, the
Cohen real $ r_{\zeta_\alpha}$ is generic
over $N_\alpha$. For this it suffices that
the countable model $N_\alpha \in V^{\P_{\zeta_\alpha}}$, which means that
$\zeta_\alpha$ just has to be sufficiently large.
Let the $a_k$ be chosen as in the proof of Lemma~4.1.
Then there are infinitely many $k$ such that
$$r_{\zeta_\alpha}^{-1}(\{\ell_{\zeta_\alpha}\}) \cap [a_{2k+\ell -1}, a_{2k+\ell}) = \emptyset,$$
and as in the proof of Lemma~4.1 this suffices.
\end{proof}

\begin{choice}\label{delphi}
We start with $\bar{M}$ as described.
By Lemma~\ref{iterate}, all the $\P_\delta$, $\delta \leq \aleph_2$, have the $\bar{M}$-\cc{}
as soon as we can arrange that all the $\Q_\delta$ have the $(\bar{M})^*$-\cc{}
in $V^{\P_\delta}$. The Cohen forcing has the $\bar{M}$-\cc{} for any $\bar{ M }$.
The  $\Q_\delta$  in the steps $\delta \in S_1^2$
can be chosen by the previous lemma so that they have the $(\bar{M})^*$-c.c.
\end{choice}

\begin{lem}\label{core}
If $\delta\in S_1^2$, $\Q_\delta$ is chosen as in Lemma~\ref{5.5},
and $\P_\delta$ satisfies (f) of Definition~\ref{5.2}, then
$\P_{\delta+1}$ has the property stated in item (f).
\end{lem}

\begin{proof}

\smallskip

Suppose that $p \Vdash_{\P_{\delta+1}}\mbox{``} \name{A} \in \roth\mbox{ and }
|\{\zeta \in W_1 \cap \delta: \name{A} \as \name{Y_\zeta^{\ell_\zeta}}\}|=
\aleph_1$'', and
w.l.o.g.\
$p \Vdash_{\P_{\delta+1}} \mbox{``} \name{A} \in \roth \mbox{ and }
\{\zeta \in W_1 \cap \delta: \name{A} \as \name{Y}_\zeta^{\ell_\zeta}\}$
is increasingly enumerated by
$ \{\xi_\alpha : \alpha < \aleph_1\}=W_1(A)$''.

\smallskip

We take for $n \in \omega$ a maximal antichain $\{p_{n,i} : i \in \omega\}$
above $p$
deciding the statements $\check{n} \in \name{A}$ with truth value $t_{n,i}$.
Let $C_{n,i} = \{ \eps \leq \delta : p_{n,i} (\eps) \neq 1\}$.
For $\eps \in C_{n,i}\cap S_1^2$ with ${\mathbb Q}_\eps \neq \{0\}$,
let $p_{n,i}(\eps) = (m_{n,i}(\eps), h_{n,i}(\eps),
F_{n,i}(\eps))$. Let $ F'_{n,i}(\eps) = \{\zeta_\alpha(\eps): \alpha \in F_{n,i}(\eps)\}$.
 We assume that all these are objects not just names.
 For $\eps \in C_{n,i}\setminus S_1^2$ let $p_{n,i}(\eps)= h_{n,i}(\eps)$,
$m_{n,i}(\eps) = |h_{n,i}(\eps)|$
and set the other two components for simplicity zero.
Set $m_{n,i} = \max\{ m_{n,i}(\eps) : \eps \in C_{n,i}\}$.
Set
\begin{multline*}
\bar{C}= \< \< (m_{n,i}(\eps),h_{n,i}(\eps), F_{n,i}(\eps), F'_{n,i}(\eps),
\< g_\alpha(\eps) \restriction m_{n,i} : \alpha \in F_{n,i}(\eps) \>):
\\
\eps \in C_{n,i} \> : n,i\in \omega\>.
\end{multline*}

For each $\beta\in \aleph_1$,
let $p_\beta \geq p$, $p_\beta \Vdash_{\P_{\delta+1}} \mbox{``}
\name{A} \cap [s_\beta,\infty) \subseteq
\name{Y}_{\xi_\beta}^{\ell_{\xi_\beta}}$''
and $p_\beta$ shall decide the value of $\ell_{\xi_\beta} \in 2$ and $s_\beta\in \omega$.
For $\beta<\aleph_1$ we set $C_\beta = \{ \eps \leq \delta : p_\beta(\eps) \neq 1\}$.
If $\eps \in C_\beta \cap S_1^2$, then $p_\beta(\eps) = (m_{\beta}(\eps), h_\beta(\eps),
F_\beta(\eps))$.
If $\eps \in C_\beta\setminus S_1^2$, then $p_\beta (\eps) = h_\beta(\eps)$,
$_\beta(\eps) =|h_\beta(\eps)|$ and $F_\beta(\eps)=\emptyset$.
For all $\beta$, $\eps\in C_\beta$, let
Let $F_\beta'(\eps) = \{\zeta_\alpha(\eps) : \alpha \in F_\beta(\eps)\}\subseteq W_2$.

Set
\begin{multline*}
R_\beta(m)= \< (m_{\beta}(\eps),h_{\beta}(\eps), F_{\beta}(\eps), F'_{\beta}(\eps),
\< g_\alpha(\eps) \restriction m : \alpha \in F_{\beta}(\eps) \>
)\\
 : \eps \in C_{\beta}\>.
\end{multline*}
These are finite arrays of finite sets.

\smallskip

Now we thin out:  First we assume that for some $k\in\omega$ for all $\beta<\aleph_1$,
$ |C_\beta|= k$, $s_\beta\leq k$.
We apply the delta system lemma to $C_\beta$, $\beta \in \aleph_1$, get a root $C$.
We assume that $\delta \in C$, as this is the difficult case.
 We apply the delta lemma for each $\eps\in C$ to the
$F_\beta(\eps)$, $\beta \in \aleph_1$, and get a root $F(\eps)$, and  to
$F'_\beta(\eps)$, $\beta \in \aleph_1$, and get a root $F'(\eps)$.
We further assume that for each $\beta$ in  the delta system
and for all $\eps \in C$,
all $F_{\beta}(\eps) \setminus F(\eps)$ are above $\max(\bigcup_{\eps'\in C}(F(\eps'))
\cup (C\setminus\{\delta\}))$
and same for the primed ones.
\comment{All $F_{\alpha_i}'(\eps)\setminus F$ are above max $F'(\eps)$.
This goes only $\eps$-wise, because in the definition of ${\mathcal K}_\gamma$ in item (e)
we did not require coherence in the enumerations
$\<\zeta_\alpha (\eps): \alpha \in \aleph_1\>$.
}
We thin out further and assume that there are $(m(\eps) ,h(\eps), F(\eps))$
such that  for all $\beta<\aleph_1$, for all $\eps \in C$,
 $m_{\beta}(\eps) = m(\eps)$,  $h_{\beta}(\eps) =
h(\eps) \in {}^{m(\eps)} \omega$, and for the $\eps \in C_\beta \setminus C$,
the increasingly enumerated $\eps$'s in $C_\beta =\{ \eps^\beta_i: i<k\}$,
are isomorphic to the lexicographically
 first $\< \eps_i : i < k \>$, i.e., $m_{\beta}(\eps^\beta_i)
= m(\eps_i)$,  $h_{\beta}(\eps^\beta_i) =
h(\eps_i) \in {}^{m(\eps_i)} \omega$,  and we use a delta system argument
 on the $F_\beta(\eps^\beta_i)$ giving a root
$F(\eps_i)$ and again impose  on the parts
$F_\beta(\eps^\beta_i)\setminus F(\eps_i)$, that they have to lie
 above $\bigcup_{i<k} F(\eps_i)$ and are all of the same size.
The analogous thinning out is done for the
primed parts, that have to lie above $\max(\bigcup_{i<k} (F'(\eps_i)) \cup (C\setminus\{\delta\}))$,
be for all $i$ of the same
size $|F'_\beta(\eps^\beta_i)|$ independently of $\beta$ (but depending on $i$),
and all of the $\< F'_\beta(\eps^\beta_i) : i<k\>$
shall have the same $\le$ or $\ge$-relations with the members of
 $C_\beta(\eps_i)$.
Moreover, if $\eps$ is a Cohen coordinate in $C_{\beta}$, then $p_{\beta} (\eps)$ does not
depend on $\beta$.

\smallskip

We let $m_{max}$ be the the maximum of the  $m(\eps)$ and of
the lengths of all the finitely many Cohen coordinates for all $\beta$ in the delta system.
Let $\triangleleft$ denote the initial segment relation for finite sequences.
We thin out further and assume that all the $R_\beta(m_{max})$
have the same quantifier free $(<_{\aleph_1}, \triangleleft)$-type over
${\rm Ran}(\bar{C})\cup {\rm Ran}({\rm Ran}(\bar{C}))$. Speaking about components of
five tuples $(m,h,F,F',\bar{ g})$ separately is allowed as well as evaluating
$\bar{g}$ and the members of all involved finite sets.
There are only countably many quantifier types in this language
that can be fulfilled by a (finite) sequence $R_\beta(m_{max})$ in our delta system.

\medskip

Let $G_{\delta}$ be a  subset of $\P_\delta$ that is generic over $V$
such  $W^* = \{\gamma \in W_1(A) \cap \delta : p_\gamma \restriction \delta \in
G_{\delta} \}$
is uncountable.

For $\gamma \in W^*$, let in $V[G_{\delta}]$,
$$B_\gamma = \{ n \in\omega
:  \exists p' \in \P_{\delta+1}, p' \geq p_\gamma,
p' \restriction \delta \in G_{\delta},\mbox{ and }
p'\Vdash_{\P_{\delta+1}} n \in A\}.$$

$B_\gamma \as Y_{\xi_\alpha}^{\ell_{\xi_\alpha}}[G]$, and the latter
 is fully evaluated by $G$,
because $\xi_\alpha \in W_1\subseteq \delta +1$ for $\alpha<\aleph_1$, and $\delta \not\in W_1$.

We shall show that for $\beta$, $\gamma \in W^*$, $B_\beta \cap [k,\infty)
 = B_\gamma \cap [k,\infty)= B\in V[G]$.
Then $B$ is a counterexample to
$\< (\P_\eps, \Q_\beta, M^\eps, W_1, W_2) : \eps \leq
\delta, \beta<\delta \> \in {\mathcal K}_\delta$.

Let $||_{\P_{\delta+1}}$ denote the compatibility relation in $\P_{\delta+1}$.
If $n \in B_\beta$, then $p_\beta ||_{\P_{\delta+1}}\, p_{n,i}$ for the one $i$
such that $p_{n,i}\in G$, and for this $i$ we have $t_{n,i}= true$.
The same holds  for $n\not\in B_\beta$ with $false$.
So our claim that $B_\beta \cap [k,\infty) =B_\gamma
\cap [k,\infty)$ for all $\beta,\gamma \in W^*$ now follows from

\begin{claim}
For all $\beta,\gamma$ in $W^*$:
$$ p_\beta ||_{\P_{\delta+1}}\, p_{n,i} \mbox{ iff } p_\gamma ||_{\P_{\delta+1}}\, p_{n,i}.$$
\end{claim}

\begin{proof} The point is the coordinate $\delta$, since the
restrictions to $\delta$ are in $G_{\delta}$, and hence compatible.
Assume $p_{n,i}(\delta) = (m_{n,i}, h_{n,i}, F_{n,i})$,
$p_\beta(\delta) = (m_\beta, h_\beta, F_\beta)$, $p_\gamma(\delta)
=(m_\gamma, h_\gamma, F_\gamma)$. We do not write the $\delta$ at these points, but
will not suppress it completely.
 We assume that $p_\beta (\delta)$ is compatible with $p_{n,i}(\delta)$.

\medskip

First case: $m_\beta \geq m_{n,i}$.
Then $p_\beta || p_{n,i}$ means
$h_\beta \triangleright h_{n,i}$  and for all $\alpha \in F_\beta \cup F_{n,i}$
for all $m\in [m_{n,i},m_\beta) \cap Y_{\zeta_\alpha(\delta)}^{\ell_{\zeta_\alpha(\delta)}}$,
$(h_\beta(m) \geq g_\alpha(\delta)(m))$.

We have to show that the same holds for $p_\gamma$.
First, by our thinning out $m_\beta =m_\gamma$,
$h_\beta = h_\gamma$, and hence $h_\gamma \triangleright h_{n,i}$,
and $F_\beta \cap F_{n,i} = F_\gamma \cap F_{n,i}$ .

1 a) We have to show: For all  $\alpha \in F_{n,i}$
for all $m\in [m_{n,i},m_\gamma) \cap Y_{\zeta_\alpha(\delta)}^{\ell_{\zeta_\alpha(\delta)}}$\\
$(h_\gamma(m) \geq g_\alpha(\delta)(m))$.

And since $h_\beta=h_\gamma$,
for all  $\alpha \in F_{n,i}$
for all $m\in [m_{n,i},m_\gamma) \cap Y_{\zeta_\alpha(\delta)}^{\ell_{\zeta_\alpha(\delta)}}$,
$(h_\gamma(m) \geq g_\alpha(\delta)(m))$.

\smallskip

1 b) We also have to show: For  all  $\alpha \in F_\gamma$
for all $m\in [m_{n,i},m_\gamma) \cap Y_{\zeta_\alpha(\delta)}^{\ell_{\zeta_\alpha(\delta)}}
(h_\gamma(m) \geq g_\alpha(\delta)(m))$.
For  $\alpha \in F_\gamma \cap F_\beta$ the latter requirement  is clearly
fulfilled, as $h_\beta = h_\gamma$.
For the part  $F_\gamma \setminus F(\delta)$ we need to look closer:
Suppose some condition in $p_\gamma$ forced something about
$Y_{\zeta_\alpha(\delta)}^{\ell_{\zeta_\alpha(\delta)}}$. Then  $p_\gamma(\zeta_\alpha(\delta))
\neq 1$ and hence $\zeta_\alpha(\delta) \in C_\gamma\cap W_2$.
But then because of the indiscernibility over $m_\gamma=m_\beta\leq m_{max}$
(which is  a component of $\bar{C}$),  $\zeta_\alpha(\delta) \in C_\beta$
and hence it is in the root $C$.
So $p_\beta$ forced by our thinning out same fact about
$ Y_{\zeta_\alpha(\delta)}^{\ell_{\zeta_\alpha(\delta)}}\cap m_{max}$.
Hence, for all $\alpha \in F_\gamma$
for all $m\in [m_{n,i},m_\gamma) \cap Y_{\zeta_\alpha(\delta)}^{\ell_{\zeta_\alpha(\delta)}}$,
$
(h_\gamma(m) \geq g_\alpha(\delta)(m))$.
So, taking 1 a) and 1 b) together,  $p_\gamma || p_{n,i}$.

\medskip

Second  case: $m_\beta \leq m_{n,i}$.
Then $h_\beta \triangleleft h_{n,i}$, and $p_\beta || p_{n,i}$ means that
for all $\alpha \in F_\beta \cup F_{n,i}$
for all $m\in [m_{\beta},m_{n,i}) \cap Y_{\zeta_\alpha(\delta)}^{\ell_{\zeta_\alpha(\delta)}}$,
$
(h_{n,i}(m) \geq g_\alpha(\delta)(m))$.
This latter statement does hold also for
$F_\gamma$ instead of  $F_\beta$ and $m_\gamma$ instead of $m_\beta$,
beause $m_\gamma=m_\beta$ and
$(F_\beta, \< g_\alpha(\delta) \restriction m_{n,i}: \alpha \in F_\beta
\>)$  and
$(F_\gamma, \< g_\alpha(\delta) \restriction m_{n,i} : \alpha \in F_\gamma
\>)$  are part of $R_{\beta}(m_{max})$ and $R_\gamma(m_{max})$ and hence
indiscernible over $h_{n,i}$ for arguments $m \in
Y_{\zeta_\alpha(\delta)}^{\ell_{\zeta_\alpha}(\delta)}$,
as for these $m$'s, that are forced to be in a Cohen part, $\zeta_\alpha(\delta) \in C$ and
hence by our thinning out we have $ m_{max} \geq m$.
Also $h_\gamma \triangleleft h_{n,i}$, and hence $p_\gamma || p_{n,i}$.

So the claim is proved and with it also Lemma~\ref{core}. \end{proof}

\begin{lem}\label{5.9}
\begin{myrules}
\item[(1)] If $\cf(\gamma) = \aleph_1$ and $\name\Q$ and $\bar{M}^{\gamma}$ are
as in the previous lemma and
if $\<  \P_\beta,
\name{\Q_\beta}, \bar{M}^\beta, W_1, W_2) :
\beta < \gamma  \>\in {\mathcal K}_\gamma$, then
$$\<\P_\beta,
 \name{\Q_\beta}, \bar{M}^\beta, W_1, W_2) :
\beta <\gamma  \>
\; \hat{} \; \< \P_\gamma, \name{\Q},
\bar{M}^{\gamma} \> \in {\mathcal K}_{\gamma+1}.$$

\item[(2)]
If $\cf(\gamma) = \aleph_0$ and
if $\<  \P_\delta,
\name{\Q_\beta}, \bar{M}^\beta, W_1, W_2) :
\beta < \gamma  \> \in {\mathcal K}_\gamma$, then
$$\<\P_\beta,
 \name{\Q_\beta}, \bar{M}^\beta, W_1, W_2) :
\beta <\gamma  \>
\; \hat{} \; \< \P_\gamma, \bbC,
\bar{M}^{\gamma} \> \in {\mathcal K}_{\gamma+1}.
$$

\item[(3)] If $\cf(\gamma) = \aleph_0$ and
if $\<  \P_\beta,
\name{\Q_\beta}, \bar{M}^\beta, W_1, W_2) :
\beta < \gamma  \> \restriction \beta
\in {\mathcal K}_\beta$ for each $\beta <\gamma$, then
$
\<\P_\beta,
 \name{\Q_\beta}, \bar{M}^\beta, W_1, W_2) :
\beta <\gamma  \> \in {\mathcal K}_\gamma$.

\item[(4)] If $\cf(\gamma) = \aleph_1$ or $\gamma=\aleph_2$, and
if $
\<  \P_\beta,
\name{\Q_\beta}, \bar{M}^\beta, W_1, W_2) :
\beta < \gamma  \> \restriction \beta
\in {\mathcal K}_\beta$ for each $\beta <\gamma$, then
$
\<\P_\beta,
 \name{\Q_\beta}, \bar{M}^\beta, W_1, W_2) :
\beta <\gamma  \>
\in {\mathcal K}_\gamma$.
\end{myrules}
\end{lem}
\begin{proof}
(1) This was proved in Lemma~\ref{core}.

(2) If $A$ is an almost subset of uncountably many $Y_\zeta$'s, then there is some $\gamma_0
<\gamma$ that there are uncountably many such $\zeta$ below $\gamma_0$.
 $A$ is possibly a name using the last, new forcing.
But this is just Cohen forcing. So there is some finite part of a Cohen condition
 forcing that $\name{A}$ is in uncountably many $Y_\zeta$'s. But then
also the forcing $\P_{\gamma}$
already contains a name for some infinite $B \subseteq \omega$
almost contained in the intersection of uncountably many $Y_\zeta$'s with $\zeta< \gamma_0$.
So $P_{\gamma}$ does not fulfill property (f) and hence the
induction hypothesis is not fulfilled.

(3) First we use the  pigeonhole principle for the $Y_\zeta$'s as in the previous item.
Then we use the following

\begin{lem}\label{5.10} Assume
\begin{myrules}
\item[(a)]  $\< {\P}_n : n \in \omega \> $ is a
$ \lessdot $-increasing sequence
of c.c.c.\ forcing notions with union ${\P}$,
\item[(b)] ${\mathcal Y} $ is a set of ${\P}_0 $-names of infinite
subsets of $\omega $,
\item[(c)]  for $ n \in \omega $ we have
$ \Vdash{}_{{\P}_n} \mbox{``} \kappa = \cf( \kappa )
> |\{ \name{ Y } \in {\mathcal Y} : \name{B}\as\name{ Y }\}|$'',
whenever $ \name{ B } $ is a $ {\P}_n $-
name of an infinite subset of $\omega $.
\end{myrules}
Then condition (c) holds for $\P$  too.
\end{lem}
\begin{proof}
Since $\P$ is a c.c.c.\ forcing notion,
also in $V^{\P}$  we have
$ \kappa $ is a regular cardinal.

If the desired conclusion fails, then we can find
a $\P$-name $ \name{ B } $ of an infinite subset of
$ \omega $  and a sequence
$ \<(p_\alpha , \name{Y}_\alpha, m_\alpha  )  : \alpha < \kappa
\> $ such that
\begin{myrules}
\item[($\alpha$)]  $m_\alpha \in \omega $,
\item[($\beta$)] $ \name{Y}_\alpha \in {\mathcal Y}  $ without repetitions,
\item[($\gamma$)] $p_\alpha \in \P$, $ p_ \alpha \Vdash_{\P}
\name{ B } \setminus m_\alpha  \subseteq \name{ Y } _ \alpha $.
\end{myrules}

Since $\cf(\kappa) > \aleph_0$, for some $ n(* ), m(*)  \in \omega $  the set
$ S =^{\rm df } \{ \alpha < \kappa :
p_ \alpha \in  {\P} _ { n(*) } , m_\alpha = m(*) \} $
has cardinality $ \kappa $. We identify it with $\kappa$.

Now for every large enough  $ \alpha \in S $ we have
$$ p_ \alpha \Vdash_ {\Bbb   P}  \kappa = |
\{ \beta \in S: p_ \beta \in
 \name{ G}_ { {\P} _ {n(*)} } \} |. $$

Why? Else for an end segment of  $\alpha <\kappa$ there is $q_\alpha \geq p_\alpha$ such that
for all but $< \kappa$ many $\beta \in S$,  $q_\alpha \Vdash p_\beta \not \in \name{G}_{\P_{n(*)}}.$
That means that for an end segments of $\alpha <\kappa$, w.l.o.g., for all
$\alpha \in \kappa$,
${\rm Perp}_\alpha := \{\beta \in S :  q_\beta \perp q_\alpha \}$
contains an end segment of $S$.
Then we take the diagonal intersection $D$ of all these end segments
of $S$. Since $\kappa$
is regular, $D$  contains a club in $\kappa$. But then $\{q_\beta : \beta \in D\}$
 is an antichain in $\P_{n(*)}$ of size $\kappa$.
Contradiction.

Let $G_{n(*)}$ be a subset of $\P_{n(*)}$
generic over $V$, and let $S_* :=
\{\beta \in S : p_\beta \in G_{n(*)}\} $. We choose $G_{n(*)}$, such that
$|S_*|= \kappa$.
We let $ B' = \cap \{\name{ Y }_\beta  \setminus m(*) :
\beta \in S_*\}$.
Then in $V[G_{n(*)}]$,
$B'$ is an infinite subset of $\omega$ included
in $\kappa$  members of ${\mathcal Y}$,
contradicting the assumption. So Lemma~\ref{5.10} is proved.
\end{proof}

(4) If $\P_\delta$ adds some $A$, then this already comes earlier, say in $V^{\P_\eps}$,
$\eps <\delta$, because $A \subseteq\omega$
and because of  the c.c.c. If $A \as Y_\zeta$ is forced,
then $\zeta <\eps$. This contradicts the induction hypothesis for $\P_\eps$.
This completes the proof of Lemma~\ref{5.9}.
\end{proof}

The lemmas together give that there is an $\aleph_2$-approximation, and the
proof of Theorem~\ref{goodclass} is completed.
\end{proof}

\medskip

With some extra care our proof can be modified to
yield the following (cf.\ \cite{ShSt465, BurMil}).

\begin{thm}\label{stronger}
It is consistent (relative to ZFC) that all of the following assertions hold:
\be
\i Each unbounded set of $\ww$ contains an unbounded subset of size $\aleph_1$,
\i Each nonmeager subset of $\ww$ contains a nonmeager subset of size
$\aleph_1$,
\i $\g =\aleph_1$; and
\i $\cov(\Dfin) = \cov(\M) = \c = \aleph_2$.
\ee
\end{thm}
\begin{proof}
This time we work with a version of ${\mathcal K}_\gamma$
with increasing oracles, which
means that the $\bar{M}^\eps$-\cc{} implies  $\bar{M}^\delta$-\cc{} for
$\eps >\delta$ and that
$\P_\delta \Vdash \mbox{``}\P_{[\delta,\eps)}\mbox{ has the }\name{\bar{M}}^{\delta+1}\mbox{-c.c.''}$,
though the initial segment need not yet fulfill it, and the name for this new oracle may not yet
have an evaluation in an initial segment $\P_\gamma$, $\gamma<\delta$.
 The new parts of the oracles take care of the
unbounded  and the nonmeager
families that appear later in the iteration and that are frozen by the next
step if their intersection with $V^{\P_\delta}$ is guessed by the diamond
sequence and happens to be unbounded or nonmeager at the current stage $\delta$:
The conservation of the unboundedness and nonmeagerness of the intersection is
written into all the oracles from $\delta$ onwards.
\end{proof}

\ed
\begin{thebibliography}{00}
\bibitem{BlSh}
A.\ R.\ Blass and S.\ Shelah,
\emph{There may be simple $P_{\aleph_1}$- and $P_{\aleph_2}$-points,
and the Rudin-Keisler ordering may be downward directed},
Annals of Pure and Applied Logic \textbf{33} (1987),
213--243.

\bibitem{BlassHBK}
A.\ R.\ Blass,
\emph{Combinatorial cardinal characteristics of the continuum},
in: \textbf{Handbook of Set Theory} (M.\ Foreman, A.\ Kanamori, and M.\ Magidor, eds.),
Kluwer Academic Publishers, Dordrecht, to appear.

\bibitem{Burke}
M.\ Burke,
\emph{Liftings for Lebesgue measure},
Israel Mathematical Conference Proceedings \textbf{6} (1993),
119--150.

\bibitem{BurMil}
M.\ Burke and A.\ W.\ Miller,
\emph{Models in which every nonmeager set is nonmeager in a nowhere dense Cantor set},
Canadian Journal of Mathematics 57 (2005), 1139--1154.
\arx{math.LO/0311443}

\bibitem{Mildenberger}
H.\ Mildenberger,
\emph{Groupwise dense families},
Archive for Mathematical Logic \textbf{40} (2001),
93--112.

\bibitem{Sh:f}
S.\ Shelah,
\emph{Proper and Improper Forcing} (second edition),
Springer, 1998.

\bibitem{ShSt465}
S.\ Shelah and J.\ Stepr\=ans,
\emph{Maximal Chains in ${}^\omega\omega$ and Ultrapowers of the Integers},
Archive for Mathematical Logic \textbf{32} (1993),
305--319.

\bibitem{ShTb768}
S.\ Shelah and B.\ Tsaban,
\emph{Critical cardinalities and additivity properties of
combinatorial notions of smallness},
Journal of Applied Analysis \textbf{9} (2003),
149--162.

\bibitem{Vaughan}
J.\ Vaughan,
\emph{Small uncountable cardinals and topology}, 
in: \textbf{Open Problems in Topology} 
(eds.\ J.\ van Mill and G.\ M.\ Reed),
North-Holland, Amsterdam: 1990,
195--218.
\end{thebibliography}
